\documentclass[12pt]{article}
\usepackage[utf8]{inputenc}   % LaTeX, comprends les accents !
\usepackage{amsmath} 
\usepackage{amssymb}
\usepackage{amsthm}
\usepackage{hyperref}

\usepackage{xcolor}
\hypersetup{
    colorlinks,
    linkcolor={red!50!black},
    citecolor={blue!50!black},
    urlcolor={blue!80!black}
}

\usepackage{tikz}
\usetikzlibrary{matrix,decorations,arrows}
\usetikzlibrary{decorations.markings}
\newcounter{remark}
\newcommand{\remark}{\addtocounter{remark}{1}
                       \par \quad {\bf \arabic{remark}}.\,
                      }

\newenvironment{rk}{\begin{quote}
                     {{\bf Remark} --}
                    }{\end{quote}}

\newenvironment{preuve}{\medbreak \noindent {\bf Proof~---}}
                       {\hfill $\square$ \medbreak}

\newcommand{\NN}{\mathbb N}
\newcommand{\ZZ}{\mathbb Z}
\newcommand{\PP}{\mathbb P}

\newcommand{\FF}{\mathbb F}
\newcommand{\OO}{\mathcal O}

\newtheorem{theo}{Theorem}%%[section]

\newtheorem{defi}[theo]{Definition}

\newtheorem{lem}[theo]{Lemma}

\newtheorem{prop}[theo]{Proposition}

\newcommand{\Aut}{\operatorname{Aut}}

\newcommand{\different}{\operatorname{\mathfrak D}}

\newcommand{\divisor}{\operatorname{div}}

\newcommand{\Id}{\operatorname{Id}}

\newcommand{\Jac}{\operatorname{Jac}}

\newcommand{\Proj}{\operatorname{Proj}}

\newenvironment{psmallmatrix}{\left(\begin{smallmatrix}}{\end{smallmatrix}\right)}

\newcommand{\Czero}{X_0}
\newcommand{\Cun}{X}
\newcommand{\Cdeux}{\Gamma}
\newcommand{\Graphe}{{\mathcal G}}

\newcommand{\ga}{\gamma}

\newcommand{\vertiii}[1]{
    {\left\vert\kern-0.25ex\left\vert\kern-0.25ex\left\vert #1 
    \right\vert\kern-0.25ex\right\vert\kern-0.25ex\right\vert}}

\newcommand{\legendre}[2]{\genfrac{(}{)}{}{}{#1}{#2}}
\newcommand{\hypergeom}[1]{F\left(\left.\genfrac{}{}{0pt}{}{1/3,\, 2/3}{1}\right|#1\right)}
\title{\bfseries A Graph Aided Strategy to Produce Good Recursive Towers over Finite Fields}
\author{
Emmanuel Hallouin \& Marc Perret\thanks{Institut de Math\'ematiques de Toulouse, UMR 5219}
}
\begin{document}
\maketitle

\begin{abstract}
We propose a systematic method to produce potentially good recursive towers over finite fields. The graph point of view, so as some
{\tt magma} and {\tt sage} computations are used in this process. We also establish some theoretical functional criterion ensuring the existence of many rational points on a recursive tower. Both points are illustrated on an example, from the production process, to the theoretical study, using this functional
criterion, of the parameters of the obtained potentially good tower.
\end{abstract}

\tableofcontents

%---------------------------
\section*{Introduction}
The search of {\em explicit} examples of sequences of algebraic curves over a given finite field, of genus growing to infinity and having as much as possible rational points with regard to their  genera became more and more important not only for its own, but also for several uses such as coding theory and cryptography (Garcia and Stichtenoth~\cite{StichtExpltTowers}) or for multiplication algorithms over finite fields (Ballet~\cite{Ballet}).
For quite a long time, only modular examples were known for square size of the finite field (Tsfasmann-Vladut-Zink~\cite{TVZ} and Ihara~\cite[``The general case'', p. 723]{Ihara}), and only examples coming from class field theory were known for non square size. Unfortunately, these examples were not explicit. 

In 1995 appeared the important Garcia-Stichtenoth's paper~\cite{GSInvent} in which the very first explicit example was given.
The explicitness comes from the recursive definition of each floor of the
tower; such towers are now called {\em recursive towers}.
 Since then, several authors gave many examples of recursive towers (see Garcia and Stichtenoth's survey~\cite{StichtExpltTowers} or Li's one~\cite{Li}). Two features appear at a look at the literature. 
First, the authors never explain how they where able to guess their examples. 
Second, once the explicit tower is given, there is always some difficulty in its study. Either the genus sequence is hard to compute (usually in the wildly ramified case), either the existence of many rational points is hard to prove (usually in the moderately ramified case). 
For instance, in the particularly interesting example of Garcia-Sitchtenoth moderately ramified tower recursively defined by the equation $y^2=\frac{x^2+1}{2x}$, the proof of the splitting behavior is quite mysteriously 
related to some functional equation satisfied by the well known Deuring polynomial\footnote{ This is the characteristic polynomial of supersingular invariants in characteristic $p$.}.

Apart from our previous article \cite{HallouinPerretMMJ}, the present work
joins in the continuation of 
Lenstra's and Beelen's one \cite{Lenstra,Beelen}.
A kind of non-existence result is proved by Lenstra~\cite{Lenstra} over prime finite fields for a very particular type of recursive towers. The proof is technical and quite intricate, at least for the authors of the present paper. A somehow understandable point is that it relies upon some {\em functional equation} labeled as~(5) therein and whose meaning is explained in a
concluding remark. Beelen proved soon after that the characteristic function of a totally splitting set has to fulfill some functional equation in case of separable variable correspondances on the projective line~\cite[Theorem 3.2]{Beelen}. Beelen also proves, in a second part of his paper, that such a functional equation can have at most one solution for a special type of recursive tower, he called {\em of type A}. This second part has already been generalized by the
authors~\cite{HallouinPerretMMJ}.

\medskip

In Section~\ref{Backgrounds}, we recall some backgrounds about recursive towers.
Geometrically, the definition only requires a base
curve~$\Cun$ and a correspondence~$\Gamma \subset \Cun\times\Cun$.
In this article we focus on the special case of recursive
towers with {\em separable variables}\footnote{In a forthcoming paper, we
intend to prove that one can always reduce the study of a general recursive
tower to the study of a recursive tower with separable variables.},
where the correspondence is given
by two morphisms~$f, g : \Cun \to\Czero$.
One of the main tool for the
study of a recursive tower is an associated graph~$\Graphe_\infty(\Cun,\Cdeux)$. 
This graph has already been introduced by Beelen
\cite{Beelen} with a slightly different definition, and has been
extensively used by the authors \cite{HallouinPerretMMJ}. 
Many features of the tower can be directly read on the graph.
In particular, the existence of some kind of finite component, the {\em $d$-regular} ones, is related
to the existence of sets of points of the base curve splitting totally
in the tower. Following Beelen (loc. cit.), this leads
to introduce the weaker notion of {\em completeness} for subset of~$\Cun(\overline{\FF_q})$.
The importance of completeness comes from Proposition~\ref{d-reg_outside-ram}, that the search for a splitting set in a tower is carried out if one knows the ramification loci of $f$ and $g$, an easy task, and if one can find complete sets, an harder task.
% roughly speaking
%points of a complete set
%are totally split in the tower (see
%Proposition~\ref{d-reg_outside-ram} for a precise statement).
%

\medbreak

In section~\ref{CompletenessCriteria}, we push further the understanding
of the functional equation and its connection with the existence of many rational points in the tower. 
%We describe this equation in a more general
%context and we explain how it is
%related to the existence of a complete set in~$\Cun$, thus to the
%set of totally splitting points in the tower. 
 Theorem~\ref{functional_completeness_criterion} asserts essentially that a finite
set of points of the base curve splits totally in the tower if and only if the characteristic function of its support satisfies some functional equation. 
%Moreover this equation turns to be explicit as soon as
%a component of the graph, called the {\em singular component}, is finite.
Beyond the fact that we do not need to assume the base curve to be the
projective line, the contribution of Section~\ref{CompletenessCriteria}
compared to Beelen's one is at first that while he proved that the functional equation is a {\em necessary} condition for {\em completeness}, we prove here that it is in fact a necessary and {\em sufficient} conditions for {\em regularness}.
Next, that the precise form of the functional equation can be easily read
from the singular graph of the tower. The sufficiency is of great importance
for the applications given in the following of the present article.%\footnote{ Theorem 3.2 in~\cite{Beelen} also contains conditions for forward an backward completeness. These can be extended in the same way completeness is in Section~\ref{CompletenessCriteria}, but we are not able to give interesting applications.}.

\medbreak

We think that Theorem~\ref{functional_completeness_criterion} is interesting
in both theoretical and algorithmic study of recursive towers.
In section~\ref{Applications}, we give several applications
that may convince the reader.

On the theoretical side, we deduce Theorem~\ref{LenstraNous} from 
Theorem~\ref{functional_completeness_criterion}, extending one of the two main
results in~\cite{HallouinPerretMMJ}.

On the explicit side, we give two applications of the functional equation. Both use also the graph theoretic approach of recursive towers. In
Section~\ref{TourGS}, we explain  how the mysterious functional equation needed in~\cite{GS_Degre2Moderee} to prove the splitting behavior of some base
points can in fact been deduced from the initial data $y^2=\frac{x^2+1}{2x}$.

Then, in Section~\ref{X3-2n}, we describe our main application of the functional equation. It consists in a systematic method to produce some potentially good recursive towers (Section~\ref{Producing}) and in original methods to
compute the asymptotic parameters of the tower. These methods use jointly the graph approach, some {\tt magma} and {\tt sage} computations, and Theorem~\ref{functional_completeness_criterion} (Sections~\ref{GenusComput} \&~\ref{LowerBound}). All these ideas are illustrated on a concrete example. An interesting feature of this example is the following. The reader will be convinced that no modular tool is used in its production. However, the tower turns to be modular (Proposition~\ref{Sur_Fp2}). This can be seen as another experimental evidence toward Elkies' modularity conjecture\footnote{Or Elkies' fantasia, as Elkies himself stated it in~\cite{ElkiesET}.} that any tame good recursive tower over a square finite field should be modular!

%%%%%%%%%%%%%Preliminaires
\section{Preliminary notations and background} \label{Backgrounds}

We recall the definition of a recursive tower and of its principal
parameters, and we associate a graph to a recursive
tower, as in our previous work  \cite{HallouinPerretMMJ}.
Some of the results of this article are also included.

\subsection{Recursive towers}

Let $\Cun$ be a smooth projective absolutely irreducible curve defined over the finite field $\FF_q$ and let~$\Cdeux \subset \Cun \times \Cun$ be an irreducible
correspondence, without no vertical nor horizontal component. The {\em singular recursive tower} ${\mathcal T}(\Cun, \Cdeux) = (X_n)_{n \geq 1}$ is defined, for $n \geq 1$, by
$$
X_n =
\left\{
(P_1, P_2, \cdots, P_n) \in \Cun^n \mid
\text{$(P_i, P_{i+1}) \in \Cdeux$ for each~$i = 1, 2, \ldots, n-1$}
\right\}.
$$
In this article, we restrict ourselves to correspondences with {\bfseries separable variables}, that is  of the form
\begin{equation} \label{Gamma-fg}
\Cdeux
=
\Cdeux_{f, g} = \{(P, Q) \in \Cun\times \Cun ; f(P)=g(Q)\},
\end{equation}
where $f, g : \Cun\to\Cun_0$ are two degree~$d$ maps from $\Cun$ to a given
smooth base curve $\Cun_0$.

The curves~$X_n$ may be singular and by desingularization process, we obtain
 the {\em smooth recursive tower}~$
\widetilde{\mathcal T}(\Cun, \Cdeux) = (\widetilde{X}_n)_{n \geq 1}$.
The interesting parameters of the towers are:
\begin{itemize}
\item the {\em arithmetic genus} of~$X_n$,~$\ga_n = \ga(X_n)$;
\item the {\em geometric genus} of~$X_n$, i.e. the {\em genus}
of~$\widetilde{X}_n$,~$g_n=g(X_n) = g(\widetilde{X}_n)$;
\item the number of rational
points over~$\FF_{q^r}$ for~$r\geq 1$,~$N_r(\widetilde{X}_n) = \sharp\widetilde{X}_n(\FF_{q^r})$;
\item last, the limit
$$
\lambda_r({\mathcal T})
=
\lim_{n\to +\infty} \frac{N_r(\widetilde{X}_n)}{g_n}.
%\quad \quad \text{and}\quad \quad
%\beta_r({\mathcal T})
%=
%\lim_{n\to +\infty} \frac{B_r(\widetilde{X}_n)}{g_n}.
$$
\end{itemize}
For a recursive tower, this limit always exists and is thus
a non negative number. The tower is said to be {\em asymptotically good} if it is non-zero for some~$r\geq 1$. This happens
%only if the number of points sequence is maximal, while the genus
%sequence is minimal. More precisely, it is well known 
(see Lemma~1 in~\cite{HallouinPerretMMJ} for instance) if and only if there exist some~$c,c' > 0$ such that
\begin{equation}\label{formes_suites_genres_pts}
g_n = c\times d^n + o(d^n)
\qquad\text{and}\qquad
N_r(\widetilde{X}_n) = c'\times d^n + o(d^n).
\end{equation}

\subsection{Graphs and Recursive towers}

To each recursive tower, one can associate a directed graph.

\begin{defi}
Let~${\mathcal T}(\Cun,\Cdeux)$ be a recursive tower, where~$\Cdeux$
is associated, via~(\ref{Gamma-fg}), to the
functions~$f,g : \Cun\to\Czero$.
\begin{enumerate}
\item The {\bfseries geometric graph}
$\Graphe_{\infty}(\Cun, \Cdeux)$
is the graph whose vertices are the geometric points
of~$\Cun$, and for which there is an oriented edge
from $P \in \Cun(\overline{\FF_q})$ to $Q \in \Cun(\overline{\FF_q})$
if $(P, Q) \in \Cdeux(\overline{\FF_q})$, that is if $f(P)=g(Q)$.
\item For~$S \subset \Cun(\overline{\FF_q})$, the
$S$-graph~${\mathcal G}_S(\Cun,\Cdeux)$
is the subgraph of~${\mathcal G}_\infty(\Cun,\Cdeux)$ whose vertices are the
points of~$S$.
\item The $r$-th {\bfseries arithmetic graph}~${\mathcal G}_r(\Cun,\Cdeux)$
is the~$\Cun(\FF_{q^r})$-graph.
\item The {\bfseries singular graph} is the union of the weakly connected
components containing a (directed) path joining a ramified point of~$f$
to a ramified point of~$g$.
\end{enumerate}
\end{defi}

The adjectives geometric, arithmetic, or singular, qualifying the different
graphs come from the following correspondence between the points of tower
and the paths in the graphs (see \cite[Section 3.2]{HallouinPerretMMJ}).

\begin{prop}\label{PointsPaths}
Let~${\mathcal T}(\Cun,\Cdeux)$ be a recursive tower, where~$\Cdeux$
is associated, via~(\ref{Gamma-fg}), to the
functions~$f,g : \Cun\to\Czero$.
There is a one-to-one correspondence between
$$
\{\text{paths of~${\mathcal G}_\infty(\Cun,\Cdeux)$ of length~$n-1$}\}
\longleftrightarrow
X_n(\overline{\FF_q}).
$$
This correspondence:
\begin{enumerate}
\item restricts to a one-to-one correspondence
between~$X_n(\FF_{q^r})$ 
and the set of paths of length~$(n-1)$ of the r-th arithmetic
graph~${\mathcal G}_r(\Cun,\Cdeux)$;
\item is such that a path corresponds to a singular point
of~$X_n$ if and only if it joins a ramified point of~$f$ to a ramified point
of~$g$.
\end{enumerate}
\end{prop}

At any point outside the ramification loci of $f$ and $g$, the in and out degrees in the geometric
graph are equal to the common degree~$d$ of both morphisms~$f,g$,
and only points that are ramified by~$f$ or~$g$ have in or out degrees less than~$d$. The following definition is of importance in this article because of the following alternative which follows from~\eqref{formes_suites_genres_pts} and~\cite[Proposition~21]{HallouinPerretMMJ}.

\begin{defi}
A finite subset~$S \subset \Cun(\overline{\FF_q})$ is said
to be {\bfseries $d$-regular} if the $S$-graph~$\Graphe_S(\Cun,\Cdeux)$
is a $d$-regular graph that is if any vertex has in and out degrees equal
to~$d$.
\end{defi}

\begin{prop}\label{Caract_Bonne_Tour}
Let~${\mathcal T}(\Cun,\Cdeux)$ be a recursive tower, where~$\Cdeux$
is associated, via~(\ref{Gamma-fg}), to the
functions~$f,g : \Cun\to\Czero$ of degree~$d$.
Then the tower~${\mathcal T}(\Cun,\Gamma)$ is asymptotically good
over~$\FF_{q^r}$ if and only if there exists some $c>0$, such that $g_n = c \times d^n+o(d^n)$, and
\begin{enumerate}
\item either the number of points of~$\widetilde{X}_n$, defined
over~$\FF_{q^r}$, coming
from the desingularization of the singular component
of~$\Graphe_r(\Cun,\Cdeux)$ has asymptotic
behaviour of the form to~$c' \times d^n + o(d^n)$ for some~$c' > 0$;
\item \label{Existe-dReg} or there exists a finite
$d$-regular set inside the $r$-th arithmetic
graph~$\Graphe_r(\Cun,\Cdeux)$.
\end{enumerate}
\end{prop}

See loc. cit. for details, but notice that if
there exists a $d$-regular
finite subset~$S \subset \Cun(\FF_{q^r})$ for some~$r\geq 1$, then the number
of paths of length~$(n-1)$ of this component is clearly~$\sharp S\times d^n$.
By Proposition~\ref{PointsPaths}, these paths are in one-to-one correspondence
with smooth points of~$X_n$ defined over~$\FF_{q^r}$.
Therefore the number of points of~$\widetilde{X}_n$ satisfies the last condition in~Proposition~\ref{Caract_Bonne_Tour}
and the tower is asymptotically good, provided that the genus sequence
do not grows too fast as required also as a first condition in~Proposition~\ref{Caract_Bonne_Tour}. Note that geometrically, the vertices/points of~$S$
are nothing else than the totally split points in
the tower. 

\medskip

One of the main result of our previous work \cite{HallouinPerretMMJ} says that
there exists at most one finite $d$-regular set. But there may exist
other interesting finite subsets satisfying a weaker property than regularness
which turns to be useful to give characterizations of $d$-regular set.

%To begin with, we state the following Proposition, which comes from~\cite{HallouinPerretMMJ} without being explicitly stated, but that we believe enlightening.
%\textcolor{red}{C'est vraiment la $d$-regularite qui est importante. Je m'en suis apercu en faisant mes exposes. Note que (i) $\Rightarrow$ (iii) est une consequence de Moscou !}
%
%\begin{prop}
%Let~$S$ be a finite subset of~$\Cun(\overline{\FF_p})$. The following
%assertions are equivalent.
%\begin{enumerate}
%\item The graph~${\mathcal G}_S$ is $d$-regular;
%\item The graph~${\mathcal G}_S$ is $d$-regular and the set~$S$ is complete;
%\item The graph~${\mathcal G}_S$ is $d$-regular and strongly connected;
%\end{enumerate}
%\end{prop}
%\begin{preuve}
%\textcolor{red}{Bla bla bla bla....}
%\end{preuve}

\begin{defi}
Let~$\Cun$ and~$\Czero$ be two smooth, projective,
absolutely irreducible curves over~$\FF_q$ and let~$f : \Cun\to\Czero$,
$g : \Cun\to\Czero$ be two morphisms of degree~$d$.
A subset~$S$  of~$\Cun(\overline{\FF_q})$ is said to be:
\begin{enumerate}
\item {\bfseries forward complete} if~$g^{-1}(f(S)) \subset S$;
%if for~$P,Q \in \Cun(\overline{\FF_q})$, one has:
%$$
%P\in S \quad \text{and} \quad f(P) = g(Q)
%\qquad\Longrightarrow\qquad
%Q \in S;
%$$
\item {\bfseries backward complete} if~$f^{-1}(g(S)) \subset S$;
\item {\bfseries complete} if it is both backward and forward
complete.
\end{enumerate}
\end{defi}

For~$S \subset \Cun(\overline{\FF_q})$ a finite subset,
if the graph~$\Graphe_S(\Cun,\Cdeux)$ is $d$-regular, then~$S$ is complete.
The converse is false, but one easily see that the following Proposition holds.

\begin{prop} \label{d-reg_outside-ram}
Let~${\mathcal T}(\Cun,\Cdeux)$ be a recursive tower, where~$\Cdeux$
is associated, via~(\ref{Gamma-fg}), to the
morphisms~$f,g : \Cun\to\Czero$ and let~$S$ be a subset of~$\Cun(\overline{\FF_q})$. Then~$S$ is~$d$-regular
if and only if~$S$ is complete and outside the ramification loci of~$f$ and~$g$.
If this is the case, then any point of~$S$ splits totally in the tower.
\end{prop}

%%%%%%%%%%%%%%%%Completeness Criteria
\section{Completeness and regularness criteria}\label{CompletenessCriteria}
We build successively in this section three characterizations of complete and of regular sets. The first one is a basic set theoretic criterion for completeness, the second one is a divisorial criterion for completeness, and the last one is a functional equation characterizing, once a first {\em complete} finite set is known, if another given finite set is {\em $d$-regular}.

%The first characterization of completeness is a purely set theoretic one. This can be seen as a generalization of~\cite[Lemma 3.1]{Beelen}, in which the condition was stated only to be necessary, and the bee curve was the projective line.
\subsection{A divisorial criterion for completeness}

\begin{lem}[Set theoretical completeness criterion]\label{set_complete_criterion}
A subset~$S \subset \Cun(\overline{\FF_q})$ is complete if and only if
there exists~$S_0 \subset \Czero(\overline{\FF_q})$ such
that~$S = f^{-1}(S_0) = g^{-1}(S_0)$.
\end{lem}

\begin{preuve}
Suppose that~$S$ is complete. First we prove that~$f(S) = g(S)$.
Let~$P\in S$. There exists some~$Q \in \Cun(\overline{\FF_q})$ such
that~$f(P) = g(Q)$. By forward completeness, one has~$Q \in S$.
This leads to~$f(S) \subset g(S)$, and the reverse inclusion is proved the
same way using backward completeness. Put~$S_0 = f(S) = g(S)$ and let us prove
that~$S = f^{-1}(S_0) = g^{-1}(S_0)$. Of
course~$S \subset f^{-1}(S_0) = f^{-1}(f(S))$. Conversely, if~$P \in f^{-1}(S_0)$
then~$f(P) \in S_0 = g(S)$, that is~$f(P) = g(Q)$ for some~$Q \in S$. By
backward completeness, we deduce that~$P \in S$ and that~$f^{-1}(S_0) \subset S$.
The proof of the converse works also in the same way.
\end{preuve}

For the remaining characterizations, we need to introduce the following notations.
\begin{itemize}
\item For~$S \subset \Cun(\overline{\FF_q})$, we
put~$\divisor(S) = \sum_{P\in S} P$ and for~$\varphi \in \FF_q(X)$, we denote
by~$\divisor(\varphi)$ the associated divisor.

\item For~$P \in \Cun(\overline{\FF_q})$ we denote by~$e_f(P)$ the
{\bfseries ramification index} of~$P$ by~$f$.

\item For~$S_0 \subset \Czero(\overline{\FF_q})$, let
$$
\different_f(S_0) = \sum_{P \in f^{-1}(S_0)} \left(e_f(P)-1\right) P
$$
be the {\bfseries restricted different} divisor of~$f$.
\end{itemize}
Under these notations, if~$S_0$ is a subset of~$\Czero(\overline{\FF_q})$,
one can point out the useful divisorial equalities
\begin{equation}\label{eq_f_etoile_div}
f^*\divisor(S_0) = \different_f(S_0) + \divisor\left(f^{-1}(S_0)\right)
\qquad\text{and}\qquad
g^*\divisor(S_0) = \different_g(S_0) + \divisor\left(g^{-1}(S_0)\right).
\end{equation}

We can then state the following, a generalization of Lenstra's identity.

\begin{prop}[Divisorial completeness criterion]\label{DivisorCriterion}
Let~$S_0$ be a finite subset of~$\Czero(\overline{\FF_q})$. The following assertions are equivalent.
\begin{enumerate}
\item The set~$f^{-1}(S_0)$ is complete.
\item The set~$g^{-1}(S_0)$ is complete.
\item $f^* \divisor(S_0) - g^*\divisor(S_0) = \different_f(S_0) - \different_g(S_0)$.
\end{enumerate}
\end{prop}
\begin{preuve}
This is a direct consequence of~\eqref{eq_f_etoile_div}
 together with the set theoretical completeness criterion (Lemma~\ref{set_complete_criterion}).
 \end{preuve}

\subsection{A functional criterion for regularness}

From this section, and for the rest of this article, given a curve $X$ defined over $\FF_q$ and $f, g \in \FF_q(X)$, we write $f \sim g$ if $\divisor(f)=\divisor(g)$ in~$\divisor(X)$, that is if there exists a constant
$c \in \FF_q^*$ such that $f=cg$.

\medskip

The characterizations of completeness given in Lemma~\ref{set_complete_criterion} and Proposition~\ref{DivisorCriterion} 
are not always sufficiently effective in practice.
The most fruitful characterization is the following functional regularness characterization.
%We start with two finite subsets~$S_0$ and~$T_0$ of~$\Czero(\overline{\FF_q})$.
%Let~$s,t\in\NN^*$ be such that~$s\sharp S_0 = t\sharp T_0$.
%We denote by~$a$ the order
%in~$\Jac(\Czero)$ of the degree zero
%divisor~$t\divisor(T_0) - s\divisor(S_0)$ and we
%consider~$\varphi \in \FF_q(\Czero)$ a function
%satisfying~$\divisor(\varphi) = a \left(t\divisor(T_0) - s\divisor(S_0)\right)$.
%Specializing to the case where the set~$f^{-1}(S_0)$ is complete one can
%deduce the following functional criterion for completeness.

\begin{theo}[functional regularness criterion]\label{functional_completeness_criterion}
Let $\Cdeux$ be a correspondence induced by~$f, g : \Cun \rightarrow \Cun_0$ as in~\eqref{Gamma-fg}. Let~$S = f^{-1}(S_0) = g^{-1}(S_0)$ be a complete subset
of~$\Cun(\overline{\FF_q})$. Let~$b$ be the order
of the degree zero divisor~$\different_f(S_0) - \different_g(S_0)$ in the jacobian~$\Jac(\Cun)$,
and let~$\rho \in \FF_q(\Cun)$
be a function
such that~$\divisor(\rho) = b(\different_f(S_0) - \different_g(S_0))$.
Let~$T_0$ be a finite subset
of~$\Czero(\overline{\FF_q})$, disjoint from the ramification locus of $g$. Let~$s,t\in\NN^*$
be such that~$t\sharp S_0 = s\sharp T_0$, let~$a$ be the order of the degree zero
divisor~$s\divisor(T_0) - t\divisor(S_0)$ in~$\Jac(\Czero)$,
and let~$\varphi \in \FF_q(\Czero)$
be a function
such that~$\divisor(\varphi) = a \left(s\divisor(T_0) - t\divisor(S_0)\right)$. 

Then the functional equation
\begin{equation}\label{EqFunct}
\rho^{ta} \times (\varphi \circ f)^b \sim (\varphi \circ g)^b 
\end{equation}
holds in~$\FF_q(\Cun_0)$ if and only if~$f^{-1}(T_0)$ is $d$-regular.
\end{theo}

\begin{rk}
In most instances\footnote{ We will prove in a forthcoming paper that the singular graph is {\em not} finite for the recursive tower $X_0({\mathcal P}_2^n)$ over ${\mathbb Q}(\sqrt{3})$, whose explicit recursive equation is given in~\cite[equations (47) and (48)]{ElkiesET}.}, the singular graph turns to be finite and complete. It is then fruitful to choose $S$ to be this singular graph, see Sections~\ref{TourGS} and~\ref{LowerBound}.
\end{rk}

\begin{preuve} The divisors of the functions~$\varphi\circ f$ and~$\varphi\circ g$ on~$\Cun$ are 
$$
\divisor(\varphi\circ f)
=
a
\left(
s f^*\divisor(T_0) - t f^*\divisor(S_0)
\right)
\qquad \hbox{and} \qquad
\divisor(\varphi\circ g)
=
a
\left(
s g^*\divisor(T_0) - t g^*\divisor(S_0)
\right),
$$
so that by difference
\begin{equation}\label{eq_div_phi_f_et_phi_g}
\divisor(\varphi\circ f)
-
\divisor(\varphi\circ g)
=
as \left(f^*\divisor(T_0) - g^*\divisor(T_0)\right)
-
at \left(f^*\divisor(S_0) - g^*\divisor(S_0)\right).
\end{equation}
Since~$S$ is complete,
it follows by the divisorial completeness criterion (Proposition~\ref{DivisorCriterion}) that $f^*\divisor(S_0) - g^*\divisor(S_0) =
\different_f(S_0)-\different_g(S_0)$.
Since none of the points of~$T_0$ are ramified by~$g$, we have~$\different_g(T_0) = 0$,
hence~(\ref{eq_f_etoile_div})
 reduces to~$f^*\divisor(T_0) - g^*\divisor(T_0)= \divisor\left(f^{-1}(T_0)\right)+\different_f(T_0) - \divisor\left(g^{-1}(T_0)\right)$.
Multiplying by~$b$ and taking~$\divisor(\rho^{ta})$ into account,
equation~(\ref{eq_div_phi_f_et_phi_g}) becomes
$$
b
\left[\divisor(\varphi\circ f)
-
\divisor(\varphi\circ g)
\right]
+ ta \divisor(\rho)
=
a
bs
\left[
\divisor\left(f^{-1}(T_0)\right)+\different_f(T_0)
-
\divisor\left(g^{-1}(T_0)\right)
\right].
$$
It follows that the functional equation~\eqref{EqFunct} holds if and only if
\begin{equation} \label{DivEqFunct}
\divisor\left(f^{-1}(T_0)\right)+\different_f(T_0)
=
\divisor\left(g^{-1}(T_0)\right).
\end{equation}
Suppose that the functional equation~\eqref{EqFunct}, or equivalently that~\eqref{DivEqFunct} do
holds. Since the effective divisor $\divisor\left(g^{-1}(T_0)\right)$
%have degree $\deg f \times \sharp(T_0)$. Hence
is reduced, the effective divisor
$\divisor\left(f^{-1}(T_0)\right)+\different_f(T_0)$ is also reduced. Moreover, the support of the divisor $\different_f(T_0)$ is contained in the support of $\divisor\left(f^{-1}(T_0)\right)$, hence we can deduce that $\different_f(T_0)=0$, that is $T_0$ is outside the ramification locus of $f$. Furthermore,~\eqref{DivEqFunct} becomes
$$\divisor\left(f^{-1}(T_0)\right) = \divisor\left(g^{-1}(T_0)\right),
$$
meaning that $f^{-1}(T_0)$ is complete by the set theoretical completeness criterion Lemma~\ref{set_complete_criterion}. It follows that $f^{-1}(T_0)$ is $d$-regular by Proposition~\ref{d-reg_outside-ram}. 
Conversely, suppose that that $f^{-1}(T_0)$ is $d$-regular, so that $T_0$ is disjoint from the ramification locus of $f$ and~$f^{-1}(T_0)$ is complete  by Proposition~\ref{d-reg_outside-ram}. Then $\different_f(T_0)=0$, and by the set theoretical completeness criterion, we have $\divisor\left(f^{-1}(T_0)\right)
=
\divisor\left(g^{-1}(T_0)\right)
$, so that~\eqref{DivEqFunct}, hence the functional equation~\eqref{EqFunct}, holds true.
\end{preuve}

\section{Applications}\label{Applications}

\subsection{Non existence of totally splitting sets for some recursive towers}
We deduce from the functional regularness criterion the following Theorem, extending the main result of~\cite{HallouinPerretMMJ}.

\begin{theo} \label{LenstraNous}
Let~${\mathcal T}(\Cun,\Cdeux)$ be a recursive tower, where~$\Cdeux$
is associated, via~(\ref{Gamma-fg}), to the
morphisms~$f,g : \Cun\to\Czero$ of degree~$d$.
Suppose that the tower is irreducible, and that there exists a complete set $S$ such that
$\different_f(S_0) = \different_g(S_0)$ for $S_0=f(S)=g(S)$.
Then the graph do not contains any non-empty finite $d$-regular component $T$ disjoint to $S$.
\end{theo}

\begin{rk}
A $d$-regular set $S$ satisfies $\different_f(S_0) = \sum_{P\in S} P = \different_g(S_0)$, so that this statement contains Theorem 19 of~\cite{HallouinPerretMMJ} in case of correspondences with separable variables. This Theorem also contains the case of a complete  loop at a point~$P$ (for instance for type A towers in Bellen's~\cite{Beelen} and in~\cite{Lenstra}), in which
case we have~$S = \{P\}$,~$S_0=\{f(P_0)\}$
and~$\different_f(S_0) =  (d-1)P=\different_g(S_0)$.
\end{rk}

\begin{preuve}
%\textcolor{red}{tout ca est devenu inutile : Let~$P$ be the vertex of the loop. The functions~$f$ and~$g$ have a common
%value~$P_0$ at~$P$, that is~$P_0 = f(P) = g(P)$. By hypothesis, the
%set~$\{P\}$ is a complete set. Therefore~$P_0$ must be totally ramified
%by~$f$ and~$g$ and~$f^* P_0 = g^*P_0 = d P$
%and~$\different_f(P_0) = \different_g(P_0) = (d-1)P$.}
 Suppose by contradiction that there do exist a finite $d$-regular component~$T$ and let~$T_0 = f(T) = g(T)$.
Since~$\Graphe_T(\Cun,\Cdeux)$ is assumed to be $d$-regular, each vertex of~$T$
is unramified by~$f$ and~$g$.
One can apply the functional regularness criterion
(Theorem~\ref{functional_completeness_criterion}) to $S_0$ and~$T_0$.
Since by assumption~$\different_f(S_0) - \different_g(S_0) = 0$, we have $b=1$ and one can choose~$\rho = 1$.
We denote by~$a$ the order, in~$\Jac(\Czero)$, of the zero degree
divisor~$\sharp S_0 \divisor(T_0) - \sharp T_0 \divisor(S_0)$, and we
consider~$\varphi \in \FF_q(\Czero)$ such
that~$\divisor(\varphi) = a(\sharp S_0 \divisor(T_0) - \sharp T_0 \divisor(S_0))$.
By the functional regularness criterion (Theorem~\ref{functional_completeness_criterion}), there
exists some~$c \in \overline{\FF_q}^*$ such that
$$
\varphi \circ g = c \times \varphi \circ f.
$$
For any~$Q,R\in\Cun(\overline{\FF_q})$ such that~$f(Q) = g(R)$, we deduce that
$$
\varphi \circ f(Q)
=
\varphi \circ g(R)
=
c\times \varphi \circ f(R).
$$
Therefore, if~$Q$ is any vertex of~$\Graphe_\infty(\Cun, \Cdeux)$, then $\varphi \circ f$ takes values in~$c^\ZZ \times \varphi\circ f(Q)$ on the vertices in the connected component of $Q$. But this set of values is finite since $c \in \overline{\FF_q}^*$.
Hence, the function~$\varphi \circ f$ takes only finitely many values
on each given connected component.
Since~$\varphi\circ f$ is a finite morphism, every connected components is thus finite.
Now, the graph~$\Graphe_\infty(\Cun, \Cdeux)$ being infinite and since since
there are only finitely many ramified points by $f$ or $g$, it must have
infinitely many finite $d$-regular
components by Proposition~\ref{d-reg_outside-ram}, 
a contradiction with Theorem~19 in~\cite{HallouinPerretMMJ} that under the irreducibility assumption, the graph contains at most one $d$-regular component.
\end{preuve}

\subsection{Understanding the splitting set of the optimal
tower~$y^2=\frac{x^2+1}{2x}$} \label{TourGS}

For some known Garcia-Stichtenoth recursive towers, especially for the tame one,
it is quite
difficult to prove that some set of places splits totally in the tower.
Suppose that the base curve is~$\PP^1$ (which is the case for all known explicit examples) whose Jacobian is trivial.
Once the characteristic function of the involved points on ${\mathbb P}^1$
is guessed (we explain in Subsection~\ref{X3-2n} how it can be guessed on an example),
the functional regularness criterion
Theorem~\ref{functional_completeness_criterion}
is the good tool to prove that they do split in the tower.
Let us illustrate this on the example of the moderate tower~${\mathcal T}\left(\PP^1, y^2=\frac{x^2+1}{2x}\right)$, studied 
for instance in~\cite{GS_Degre2Moderee}
and well known to be optimal. The genus computation being not difficult in this case, the hard point is the existence of some totally splitting locus.
Let~$H$ be the Deuring polynomial over $\FF_p$, whose simple roots are supersingular $j$-invariants in characteristic $p$. The splitting behaviour of some set closely related to the zero set of $H$ is easily deduced by Garcia and
Stichtenoth from the functional equation
\begin{equation} \label{EqufonctDeuring}
H(x^4)=x^{p-1}H\left(\left(\frac{x^2+1}{2x}\right)^2\right),
\end{equation}
which seems to be pulled out of the hat! We explain in this subsection that  taking into account the singular part of the graph, this is nothing but the functional equation~\eqref{EqFunct} requested in the functional regularness criterion (Theorem~\ref{functional_completeness_criterion}) for the characteristic function of some splitting set!

\medskip

\begin{figure}
\centering
\begin{tikzpicture}[line width = 1pt]
\fill (0,0) circle (3pt) node[left] {$1$} ;
\fill (2,0) circle (3pt) node[above left] {$-1$} ;
\fill (4,1) circle (3pt) node[above] {$\imath$} ;
\fill (4,-1) circle (3pt) node[below] {$-\imath$} ;
\fill (6,0) circle (3pt) node[above right] {$0$} ;
\fill (8,0) circle (3pt) node[right] {$\infty$} ;
\draw (7,1) node {$\imath^2 = -1$};
\draw[decoration={markings, mark=at position 0.625 with {\arrow{>}}},
      postaction={decorate}]
     (-0.5,0) circle (0.5cm) ;
\draw[decoration={markings, mark=at position 0.5 with {\arrow{>}}},
      postaction={decorate}] (0,0) to (2,0) ;
\draw[decoration={markings, mark=at position 0.5 with {\arrow{>}}},
      postaction={decorate}] (2,0) to [bend left] (4,1) ;
\draw[decoration={markings, mark=at position 0.5 with {\arrow{>}}},
      postaction={decorate}] (2,0) to [bend right] (4,-1) ;
\draw[decoration={markings, mark=at position 0.5 with {\arrow{>}}},
      postaction={decorate}] (4,1) to [bend left] (6,0) ;
\draw[decoration={markings, mark=at position 0.5 with {\arrow{>}}},
      postaction={decorate}] (4,-1) to [bend right] (6,0) ;
\draw[decoration={markings, mark=at position 0.5 with {\arrow{>}}},
      postaction={decorate}] (6,0) to (8,0) ;
\draw[decoration={markings, mark=at position 0.2 with {\arrow{>}}},
      postaction={decorate}]
     (8.5,0) circle (0.5cm) ;
\end{tikzpicture}
\caption{The singular component of~${\mathcal T}\left(\PP^1, \frac{x^2+1}{2x} = y^2\right)$ over~$\FF_{p^2}$
}
\label{graphe_preferee}
\end{figure}
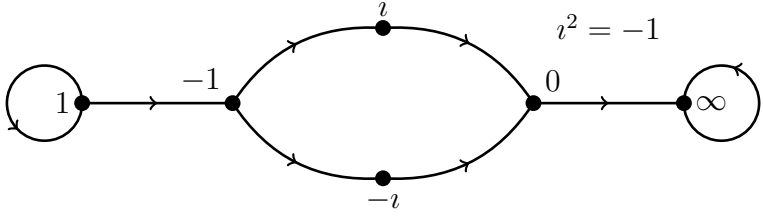
Let~$h(x) = \prod_{P\in T_0}(x-x(P)) \in \FF_q[x]$ be the characteristic polynomial of any finite set $T_0\subset \PP^1 \setminus \{\infty\}$ outside the ramification locus $\{0, \infty\}$ of~$g(y)=y^2$.

The singular graph is drawn in figure~\ref{graphe_preferee}.
For general $p\geq 3$, the involved points of $\Cun = {\mathbb P}^1$ are $S=\{1, -1, \imath, -\imath, 0, \infty\}$ where $\imath$ is a square root of $-1$ in $\overline{\FF_p}$. Applying $f(x)=\frac{1+x^2}{2x}$ (or $g(y)=y^2$), we have $S_0=\{1, -1, 0, \infty\}$. Moreover,  $\different_f(S_0)= (2-1)[1]+(2-1)[-1]$ and
$\different_g(S_0)=(2-1)[0]+(2-1)[\infty]$, hence the condition~$\divisor(\rho) = \different_f(S_0)-\different_g(S_0) = [1]+[-1]-[0]-[\infty]$ is fulfilled by the function 
$$
\rho(x) = \frac{(x-1)(x+1)}{x} \in \FF_p(\PP^1).
$$
Some {\tt magma} experiments, for few small values of~$p$, show
that there do exist such a set with~$\sharp T_0 =p-1$. Since~$\sharp S_0 = 4$ and~$p$ is odd, 
one can try~$t = \frac{p-1}{2}$ and~$s = 2$. The functional regularness criterion also requires a function~$\varphi\in \FF_p(\Czero)$
such that~$\divisor(\varphi)=2\divisor(T_0)-\frac{p-1}{2}\divisor(S_0)$.
Up to a constant, 
$$
\varphi(x)= \frac{h(x)^2}{\left[(x-1)(x+1)x\right]^{\frac{p-1}{2}}}
$$
works.
Hence, Theorem~\ref{functional_completeness_criterion} asserts that the
set $f^{-1}(T_0)$ is complete if and only if~$h(0) \not=0$ and
$$
\rho(x)^{\frac{p-1}{2}} \varphi\left(\frac{x^2+1}{2x}\right) \sim \varphi(x^2).
$$
%The condition $\varphi(\infty)\neq 0$ would holds since $h$ is where a polynomial, and $\varphi(0) \neq 0$ means that $h(0)\neq 0$.
This functional equation can be written
$$
\left[ 
\frac{(x-1)(x+1)}{x}
\right]^{\frac{p-1}{2}}
\times
\frac
{h\left(\frac{x^2+1}{2x}\right)^2}
{
\left[
\left(\frac{x^2+1}{2x}-1\right)\left(\frac{x^2+1}{2x}+1\right) \left(\frac{x^2+1}{2x}\right) 
\right]^{\frac{p-1}{2}}
}
\sim 
\frac{h(x^2)^2}
{\left[(x^2-1)(x^2+1)x^2\right]^{\frac{p-1}{2}}}
.
$$
After simplification by powers of $x$, $x-1$ and $x+1$, this is equivalent to
$$
\left[
x^{p-1}h\left(\frac{x^2+1}{2x}\right)
\right]^2
\sim
\left[
h(x^2)
\right]^2,
$$
that is to
$$
x^{p-1}h\left(\frac{x^2+1}{2x}\right)\sim h(x^2),
$$
which is neither than~\eqref{EqufonctDeuring} with $h(x)=H(x^2)$.

\medskip

Of course, proving that this functional equation do have a solution $H$ is another task. On this example, the solution~$H$ have been already guessed by Garcia and Stichtenoth. To ovoid cheating, we explain in the following section how this task can be achieved on another example, chosen in such a way that the solution is {\em a priori} unknown.

\subsection{Graph based strategy to produce and to study asymptotically good recursive towers} \label{X3-2n}
%Let us search for an explicit pair $(\Cun, \Cdeux)$ such that $\lambda(\Cun, \Cdeux)>0$ using the graph point of vue. At a first try, we choose for $\Cun$ the projective line ${\mathbb P}^1$. Thanks to proposition ?? in~\cite{HallouinPerretMMJ}, given a correspondence  $\Cdeux \sim dH+dV \in \NS({\mathbb P}^1\times {\mathbb P}^1)$ satisfying assumption (*) on p. ?? in~\cite{HallouinPerretMMJ}, the number of rational points at the $n$-th level will be at least
%$$
%\sharp \widetilde{C}_n(\FF_{q^r}) \geq \sharp S . d^n|
%$$
%if there exists a finite totally spliting  component $\Graphe_S$ in $\Graphe_r$ for some $S \subset {\mathbb P}^1(\FF_{q^r})$. Moreover, because of proposition ?? in~\cite{HallouinPerretMMJ}, it is necessary for the tower to be singular. We also need of course to be able to deal with all singular points in order to compute their contribution to the genus $g_n$. From proposition ?? in~\cite{HallouinPerretMMJ}, these singular points correspond to path in $\Graphe_s$ joining a non-etale point of $g$ to a non-etale point of $f$.  The computations will likely be easier if  there are a limited number of type of singularities.  This will be done below by targeting a particularly simple singular graph given in figure~\ref{GrapheSingVise}.

%In this section, we firstly illustrate  how the graph point of vue, together with some {\tt magma} computing help, can be used
%to generate candidates of good recursive towers. Secondly, we  illustrate how Theorem~\ref{functional_completeness_criterion} can be used to check that this candidate do works. 

The goal of this section is to show how a good recursive tower can be studied, from its production process in Section~\ref{Producing}, to the computation of its parameters in Sections~\ref{GenusComput} and~\ref{LowerBound}.

\subsubsection{Producing a potentially good recursive tower using graphs and computer help}\label{Producing}
The strategy to produce such a candidate is the following.
For algorithmic purpose, we fix a base curve~$\Cun$ and a ``kind'' of correspondence~$\Cdeux \subset \Cun\times\Cun$, that is a set of correspondences parametrized by some quasi-projective variety $V$. In order to obtain recursive towers having potentially low genus sequence as required by~Proposition~\ref{Caract_Bonne_Tour}, we also fix
a specific finite singular graph. Now, any edge of the graph ${\mathcal G}_{\infty}(\Cun, \Cdeux)$ corresponds to a relation on the parameters of the correspondence $\Cdeux$  viewed as a points of $V$. Hence, correspondences $\Cdeux$ on $\Cun$ of the given kind such that
${\mathcal G}_{\infty}(\Cun, \Cdeux)$ contains the fixed finite graph as a subgraph are parametrized by a subvariety $W$ of $V$.
Then, a {\tt magma} program returns, for a given prime~$p$, all $\FF_p$-rational points of $W$.

After very few trials of specific singular graphs\footnote{Trying with the singular graph in figure~\ref{graphe_preferee} gives as an unique solution (up to $\Aut(\PP^1)$) the Garcia-Stichtenoth tower $y^2=\frac{x^2+1}{2x}$!}, we obtain for few values of~$p$ some explicit equations. The associated graph having, as checked using a {\tt sage} program, some $d$-regular component, these equations define {\em potentially good} recursive towers. Indeed, they experimentally possess some $d$-regular set ensuring that the number of points is large enough, and their singular graph have a certain imposed shape, making it possible that the genus sequence is low enough as required in~Proposition~\ref{Caract_Bonne_Tour}.

%\subsubsection{Choosing a singular graph}

\medskip

Here is how all this works on an example. We choose:
\begin{itemize}
\item the base curve~$\Cun$ to be~$\PP^1$;
\item correspondences~$\Cdeux$ of bi-degree~$(2,2)$ and with
separable variables, that is of type~$\Cdeux_{f,g}$ where~$f$ and~$g$ are
two functions from~$\PP^1$ to~$\PP^1$ of degree~$2$;
\item a singular graph of the form%the graph of figure~\ref{GrapheSingVise}
%\begin{figure}
%\centering
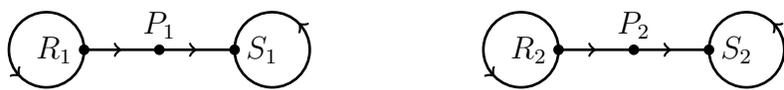
\begin{center}
\begin{tikzpicture}[line width = 1pt]
\draw[decoration={markings, mark=at position 0.625 with {\arrow{>}}},
      postaction={decorate}]
     (0,0) circle (0.5cm) ;
\draw[decoration={markings, mark=at position 0.25 with {\arrow{>}}},
      decoration={markings, mark=at position 0.75 with {\arrow{>}}},
      postaction={decorate}] (0.5,0) -- (2.5,0) ;
\draw[decoration={markings, mark=at position 0.125 with {\arrow{>}}},
      postaction={decorate}] (3,0) circle (0.5cm) ;
\fill (0.5,0) circle (2pt) node[left] {$R_1$} ;
\fill (1.5,0) circle (2pt) node[above] {$P_1$} ;
\fill (2.5,0) circle (2pt) node[right] {$S_1$} ;
\end{tikzpicture}
\hspace{2cm}
\begin{tikzpicture}[line width = 1pt,]
\draw[decoration={markings, mark=at position 0.625 with {\arrow{>}}},
      postaction={decorate}]
     (0,0) circle (0.5cm) ;
\draw[decoration={markings, mark=at position 0.25 with {\arrow{>}}},
      decoration={markings, mark=at position 0.75 with {\arrow{>}}},
      postaction={decorate}] (0.5,0) -- (2.5,0) ;
\draw[decoration={markings, mark=at position 0.125 with {\arrow{>}}},
      postaction={decorate}] (3,0) circle (0.5cm) ;
\fill (0.5,0) circle (2pt) node[left] {$R_2$} ;
\fill (1.5,0) circle (2pt) node[above] {$P_2$} ;
\fill (2.5,0) circle (2pt) node[right] {$S_2$} ;
\end{tikzpicture},
\end{center}
%\caption{The target singular graph}
%\label{GrapheSingVise}
%\end{figure}
where~$R_1,R_2$ are the ramified points of~$f$
and~$S_1,S_2$ those of~$g$. In particular, in view of the specific shape of the singular graph, these
four ramification points should be distinct.
\end{itemize}

\begin{rk}
This kind of singular graph is not chosen at random.
In order to minimize the genus sequence, the curves~$X_n$
must be singular as pointed out in~\cite{HallouinPerretMMJ}. This means by Proposition~\ref{PointsPaths} that there must exist some 
paths from ramification
points of~$f$ to ramification points of~$g$. We also have noted that in most
known examples of good recursive tower, there are loops at some of these ramification points. The chosen
graph is one of the most simple that takes into account these constraints.
\end{rk}

Note that there are $6$ parameters in such a singular graph, namely the $6$ points $P_i, R_i, S_i$, for $i = 1, 2$. We use the automorphisms group $\Aut(\PP^1)$ to fix some of them as follows. For any~$\sigma,\tau\in\Aut(\PP^1)$, the map
$$
(P_1,\ldots,P_n) \mapsto (\tau^{-1}(P_1),\ldots,\tau^{-1}(P_n))
$$
defines an isomorphism from the tower~${\mathcal T}(\PP^1, \Cdeux_{f,g})$
to the tower~${\mathcal T}(\PP^1, \Cdeux_{\sigma\circ f\circ\tau,\sigma\circ g\circ\tau})$.
Therefore, using the simply $3$-transitivity of~$\Aut(\PP^1)$
on~$\PP^1$, one can suppose that~$R_1 = 1$,
$S_1 = 0$ and~$S_2 = \infty$. In the same way, one can suppose
that~$g(S_1) = 0$, $g(S_2) = \infty$ and~$g(1) = 1$.
These normalizations lead to~$g(y) = y^2$,
so that we are reduced to look for a rational function
$$
f(x) = \frac{a_2x^2 + a_1x + a_0}{b_2x^2 + b_1x + b_0}
$$
where the parameter~$(a_2:a_1:a_0:b_2:b_1:b_0)$ lives in~$\PP^5$. Since~$f$
have to be of degree~$\deg f=\deg g = 2$, this point in~$\PP^5$ must lie in the
complementary $V$
of the sets of constant functions and of degree one functions. The constant functions is the zero set of the $2$-by-$2$
minors of~$\begin{psmallmatrix}a_2&a_1&a_0\\b_2&b_1&b_0\end{psmallmatrix}$, and the degree one functions is the zero set of the resultant
of~$a_2x^2 + a_1x + a_0$ and~$b_2x^2 + b_1x + b_0$.

For such a $f \in V$, the equation of the correspondence~$\Cdeux_{f,g}$
in~$\PP^1\times\PP^1 = \Proj(\FF_q[X_1,Y_1]) \times \Proj(\FF_q[X_2,Y_2])$
is
$$
E(X_1,Y_1,X_2,Y_2)
=
Y_2^2 (a_2X_1^2 + a_1 X_1Y_1 + a_0 Y_1^2)
-
X_2^2 (b_2X_1^2 + b_1 X_1Y_1 + b_0 Y_1^2)
=
0.
$$
%uld be related by path in
%
%
%omponent
%{GrapheSingVise}
%\Cun$.

Let us now traduce as equations on the parameter~$(a_2:a_1:a_0:b_2:b_1:b_0)$  the requested shape for the singular graph. The fact that the point~$(1:1)$ is ramified by~$f$ gives as a first
constraint the vanishing at $x=1$ of the derivative of $f(x)$, a quadratic equation in the parameters. Note for later use that the second ramified point by~$f$ is~$(a_1b_0-a_0b_1 : a_2b_1-a_1b_2)$.
Each loop, so as each horizontal path of length $2$, gives rise to linear algebraic constraints on the parameters as follows. The
four loops lead to the equations
$$
E(1,1,1,1)
=
E(0,1,0,1)
=
E(1,0,1,0)
=
E(a_1b_0-a_0b_1, a_2b_1-a_1b_2,a_1b_0-a_0b_1, a_2b_1-a_1b_2) = 0.
$$
The existence of a  path of length $2$ from $R_i$ to $S_i$ is equivalent to the vanishing of a resultant, since 
this is the request of the existence of some common zero $P_i$ of $f(R_i)-g(P)$ and of $f(P)-g(S_i)$.
Finally these $1+4+2=7$ equations define a quasi-projective curve $W$ in~$\PP^5$. 

For the small
values of the prime~$p \in \{5, 7, 11, 13, 17\}$, we find  using {\tt magma}
the unique
solution whose equation~$f(x) = \frac{x^2+x}{3x-1}$ doesn't depend on~$p$.

\medbreak

To conclude, this graph aided strategy suggests
to study the recursive tower defined by~$\Cun = \PP^1$ for $p \geq 5$, and
$$
\Cdeux_{f,g}
=
\left\{
\left((x_1:y_1),(x_2:y_2)\right) \in \PP^1\times\PP^1
\mid
\frac{x_1^2 + x_1 y_1}{3x_1 y_1 -y_1^2}
=
\frac{x_2^2}{y_2^2}
\right\}.
$$
The ramification points of the function~$f(x)=\frac{x^2+x}{3x-1}$ are~$1$
and~$-\frac{1}{3} \in \FF_p$,
those of the function~$g(y)=y^2$ are~$0$ and~$\infty$,
and the singular complete subgraph is the following 
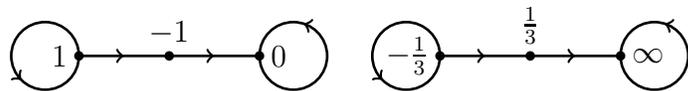
\begin{center}
\begin{tikzpicture}[line width = 1pt, scale = 0.6]
\fill (0,0) circle (3pt) node[left] {$1$} ;
\fill (2,0) circle (3pt) node[above] {$-1$} ;
\fill (4,0) circle (3pt) node[right] {$0$} ;
\draw[decoration={markings, mark=at position 0.625 with {\arrow{>}}},
      postaction={decorate}]
     (-0.75,0) circle (0.75cm) ;
\draw[decoration={markings, mark=at position 0.5 with {\arrow{>}}},
      postaction={decorate}] (0,0) to (2,0) ;
\draw[decoration={markings, mark=at position 0.5 with {\arrow{>}}},
      postaction={decorate}] (2,0) to (4,0) ;
\draw[decoration={markings, mark=at position 0.2 with {\arrow{>}}},
      postaction={decorate}]
     (4.75,0) circle (0.75cm) ;
\fill (8,0) circle (3pt) node[left] {$-\frac{1}{3}$} ;
\fill (10,0) circle (3pt) node[above] {$\frac{1}{3}$} ;
\fill (12,0) circle (3pt) node[right] {$\infty$} ;
\draw[decoration={markings, mark=at position 0.625 with {\arrow{>}}},
      postaction={decorate}]
     (7.25,0) circle (0.75cm) ;
\draw[decoration={markings, mark=at position 0.5 with {\arrow{>}}},
      postaction={decorate}] (8,0) to (10,0) ;
\draw[decoration={markings, mark=at position 0.5 with {\arrow{>}}},
      postaction={decorate}] (10,0) to (12,0) ;
\draw[decoration={markings, mark=at position 0.2 with {\arrow{>}}},
      postaction={decorate}]
     (12.75,0) circle (0.75cm) ;
\end{tikzpicture}.
\end{center}
The remaining of this section is devoted to the proof
of the following result stating that this {\em potentially good} tower is
actually {\em asymptotically good}.

\begin{theo} \label{NewTowerIsGood}
The recursive tower~$
{\mathcal T}\left(\PP^1, y^2 = \frac{x^2+x}{3x-1}\right)$
is asymptotically good.
\end{theo}

\begin{preuve}
The proof is divided in two steps. In Propositions~\ref{Genus_Sequence},
the genus sequence is computed and we observe that it has the good
shape requested in Proposition~\ref{Caract_Bonne_Tour}. Then, we 
prove~Proposition~\ref{number_of_points_sequence} that the number of
points sequence also have the good shape. 
\end{preuve}

\subsubsection{Genus sequence computation}\label{GenusComput}

In this section we establish a closed formula for the genus
sequence of the tower. One can distinguish at least two strategies to compute
the genus sequence of such a recursive tower.
The first one, used by Garcia and Stichtenoth in
their articles in this area, consists in applying the Riemann-Hurwitz
genus formula to the function field extensions~$\FF_{p}(X_n)/\FF_p(X_1)$
after having computed the different divisor of these extensions.
The second one, used in our previous work
on recursive towers~\cite{HallouinPerretMMJ},
consists in computing the geometric genus
of the curves~$X_n$ from their arithmetic one by subtracting the sum
of the measures of singularity of the points. This last point of view
 takes better into account the geometry of the data.
Furthermore, as will be seen by the reader in the proof of Lemma~\ref{Genus_Sequence} below, 
 this method is a guide for the choice of a fixed singular graph as was done in Subsection~\ref{Producing}.
This is the first time that we illustrate our method on a simple, though
non-trivial, example of recursive tower.

\begin{prop} \label{Genus_Sequence}
The genus sequence~$(g_n)_{n\geq 1}$ of the
tower~${\mathcal T}\left(\PP^1, y^2 = \frac{x^2+x}{3x-1}\right)$
is given by
$$
g_n = 2^n - (2 + n\bmod 2) \times 2^{\lfloor\frac{n}{2}\rfloor} + 1,
\qquad \forall n\geq 1.
$$
\end{prop}

\begin{preuve}
Let~$n\geq 1$. We denote by~$\nu_n$ the desingularization
morphism~$\nu_n : \widetilde{X}_n \to X_n$ and by~$X_n^\sharp$
the pullback
of the embedding~$X_n \hookrightarrow X_{n-1}\times X$
along~$\nu_{n-1} \times \Id : \widetilde{X}_{n-1}\times\Cun \to X_{n-1}\times\Cun$.
We have the cartesian diagram
$$
\begin{tikzpicture}[>=latex]
\matrix[matrix of math nodes,row sep=0.75cm,column sep=0.75cm]
{
|(Xdieze)| X_n^\sharp    & |(XtildeX)| \widetilde{X}_{n-1}\times\Cun \\
|(X)|      X_n         & |(XX)| X_{n-1}\times\Cun \\
};
\draw[->] (Xdieze) -- (X) ;
\draw[->] (XtildeX) -- (XX) node[midway,right] {\tiny $\nu_{n-1} \times \Id$};
\draw[right hook->] (Xdieze) -- (XtildeX) ;
\draw[right hook->] (X) -- (XX) ;
\end{tikzpicture}
$$
The curves~$X_n,X_n^\sharp$ and~$\widetilde{X}_n$ are birational
and the two first one are singular for~$n\geq 3$. We thus
introduce the {\em measure of singularity} of the curve~$X_n^\sharp$
defined, as usual, by
$$
\Delta_n
=
\sum_{P \in X_n^\sharp(\overline{\FF_p})} \dim_{\overline{\FF_p}} \widetilde{\OO}_P / \OO_P,
$$
where~$\OO_P$ and~$\widetilde{\OO}_P$ denote the local ring at~$P$ 
of~$X_n^\sharp$ and of its integral closure. 
Then, using Proposition~4 in~\cite{HallouinPerretMMJ}
specialized to the case where~$d=2$, $g_1 = 0$
and~$\gamma_2 = 1$, we obtain:
$$
g_n = 1 + (n-2)2^{n-1} - \sum_{i=2}^n 2^{n-i}\Delta_i.
$$
From the following Lemma~\ref{AnneauxLocaux}, we have~$\Delta_i = 2^{i-1} - 2^{\left\lfloor\frac{i}{2}\right\rfloor}$. Proposition~\ref{Genus_Sequence} follows from an easy summation exercise.
\end{preuve}

\begin{lem} \label{AnneauxLocaux}
Let~$n \geq 1$. The singular points of the curve~$X_n^\sharp$ are above the
points~$\left(1^r, -1, 0^s\right)$
and~$\left(-\frac{1}{3}^r, \frac{1}{3}, \infty^s\right)$ of~$X_n$,
for~$r+s+1=n$ and~$r \geq s \geq 1$. For such~$r,s$, there are exactly~$2^{s-1}$
points on~$X_n^\sharp$, all of them giving rise to two points
on~$\widetilde{X}_n$, and having a measure of singularity
equal to~$2^{r-s}$. The global measure of singularity of the curve~$X_n^\sharp$
is~$\Delta_n = 2^{n-1} - 2^{\left\lfloor\frac{n}{2}\right\rfloor}$.
\end{lem}

\begin{rk}
The curves~$X_1$ and~$X_2$ are smooth and thus so are~$X_1^\sharp$ and~$X_2^\sharp$.
\end{rk}

\begin{preuve}
The singular points of~$X_n$ correspond by Proposition~\ref{PointsPaths} to the
paths of the graph~${\mathcal G}_\infty$ joining a ramified point
of~$f$ to a ramified point of~$g$. In view of the chosen singular graph,
they are the points~$\left(1^r, -1, 0^s\right)$
and~$\left(-\frac{1}{3}^r, \frac{1}{3}, \infty^s\right)$
for~$r,s\geq 1$ and~$r+s+1=n$. The second one going the same way, let us concentrate on the first type of
points.
For each edge of the corresponding singular component of the singular
graph
$$
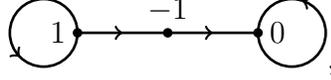
\begin{tikzpicture}[line width = 1pt, scale = 0.6]
\fill (0,0) circle (3pt) node[left] {$1$} ;
\fill (2,0) circle (3pt) node[above] {$-1$} ;
\fill (4,0) circle (3pt) node[right] {$0$} ;
\draw[decoration={markings, mark=at position 0.625 with {\arrow{>}}},
      postaction={decorate}]
     (-0.75,0) circle (0.75cm) ;
\draw[decoration={markings, mark=at position 0.5 with {\arrow{>}}},
      postaction={decorate}] (0,0) to (2,0) ;
\draw[decoration={markings, mark=at position 0.5 with {\arrow{>}}},
      postaction={decorate}] (2,0) to (4,0) ;
\draw[decoration={markings, mark=at position 0.2 with {\arrow{>}}},
      postaction={decorate}]
     (4.75,0) circle (0.75cm) ;
\end{tikzpicture},
$$
one can associate an appropriate recursive algebraic relation and
the corresponding Newton polygon in the following way:
$$
{\renewcommand{\arraystretch}{2}
\begin{array}{c|c||c|c}
\hline
&(x_n-1)^2 - 2(x_n-1) - \frac{\left(x_{n-1}-1\right)^2}{3x_{n-1}-1}
&
&(x_n+1)^2 - 2(x_n+1) - \frac{\left(x_{n-1}-1\right)^2}{3x_{n-1}-1} \\
\begin{tikzpicture}[line width = 1pt, scale = 0.6]
\fill (0,0) circle (3pt) node[left] {$1$} ;
\draw[decoration={markings, mark=at position 0.625 with {\arrow{>}}},
      postaction={decorate}]
     (-0.75,0) circle (0.75cm) ;
\end{tikzpicture}
&
\begin{tikzpicture}[>=latex]
\draw[->] (-0.5,0) -- (2.5,0);
\draw[->] (0,-0.5) -- (0,1.5);
\draw[very thick] (0,1) -- (1,0)
node[midway,anchor=south west] {$-2v(x_{n-1}-1)$};
\draw[very thick] (1,0) -- (2,0);
\draw (0,1) node {$\bullet$} node[left,anchor=east] {$2v(x_{n-1}-1)$};
\draw (1,0) node {$\bullet$} node[below] {$1$};
\draw (2,0) node {$\bullet$} node[below] {$2$};
\draw (0,0) node[below left] {$O$};
\end{tikzpicture}
&
\begin{tikzpicture}[line width = 1pt, scale = 0.6]
\fill (0,0) circle (3pt) node[above] {$1$} ;
\fill (2,0) circle (3pt) node[above] {$-1$} ;
\draw[decoration={markings, mark=at position 0.5 with {\arrow{>}}},
      postaction={decorate}] (0,0) to (2,0) ;
\end{tikzpicture}
&
\begin{tikzpicture}[>=latex]
\draw[->] (-0.5,0) -- (2.5,0);
\draw[->] (0,-0.5) -- (0,1.5);
\draw[very thick] (0,1) -- (1,0)
node[midway,anchor=south west] {$-2v(x_{n-1}-1)$};
\draw[very thick] (1,0) -- (2,0);
\draw (0,1) node {$\bullet$} node[left,anchor=east] {$2v(x_{n-1}-1)$};
\draw (1,0) node {$\bullet$} node[below] {$1$};
\draw (2,0) node {$\bullet$} node[below] {$2$};
\draw (0,0) node[below left] {$O$};
\end{tikzpicture} \\
\hline
&x_n^2 - \frac{x_{n-1}\left(x_{n-1}+1\right)}{3x_{n-1}-1}
&
&x_n^2 - \frac{x_{n-1}\left(x_{n-1}+1\right)}{3x_{n-1}-1}\\
\begin{tikzpicture}[line width = 1pt, scale = 0.6]
\fill (0,0) circle (3pt) node[above] {$-1$} ;
\fill (2,0) circle (3pt) node[above] {$0$} ;
\draw[decoration={markings, mark=at position 0.5 with {\arrow{>}}},
      postaction={decorate}] (0,0) to (2,0) ;
\end{tikzpicture}
&
\begin{tikzpicture}[>=latex]
\draw[->] (-0.5,0) -- (2.5,0);
\draw[->] (0,-0.5) -- (0,1.5);
\draw[very thick] (0,1) -- (2,0) node[midway,anchor=south west] {$-\frac{1}{2}v(x_{n-1}+1)$};
\draw (0,1) node {$\bullet$} node[left,anchor=east] {$v(x_{n-1}+1)$};
\draw (2,0) node {$\bullet$} node[below] {$2$};
\draw (0,0) node[below left] {$O$};
\end{tikzpicture}
&
\begin{tikzpicture}[line width = 1pt, scale = 0.6]
\fill (0,0) circle (3pt) node[right] {$0$} ;
\draw[decoration={markings, mark=at position 0.2 with {\arrow{>}}},
      postaction={decorate}]
     (0.75,0) circle (0.75cm) ;
\end{tikzpicture}
&
\begin{tikzpicture}[>=latex]
\draw[->] (-0.5,0) -- (2.5,0);
\draw[->] (0,-0.5) -- (0,1.5);
\draw[very thick] (0,1) -- (2,0) node[midway,anchor=south west] {$-\frac{1}{2}v(x_{n-1})$};
\draw (0,1) node {$\bullet$} node[left,anchor=east] {$v(x_{n-1})$};
\draw (2,0) node {$\bullet$} node[below] {$2$};
\draw (0,0) node[below left] {$O$};
\end{tikzpicture}\\
\hline
\end{array}
}
$$
In this tabular,~$v$ respectively denotes a valuation of the
field~$\FF_p(X_{n-1})$ satisfying~$v(x_{n-1}-1) > 0$, $v(x_{n-1}-1) > 0$,
$v(x_{n-1}+1) > 0$ and~$v(x_{n-1}) > 0$. In the remaining of this section,
for~$r,s\geq 1$, we denote by~$v_{(1^r,-1,0^s)}$ any valuation of the
field~$\FF_p(X_{r+s+1})$ satisfying~$v_{(1^r,-1,0^s)}(x_i-1) > 0$
for~$1\leq i \leq r$,
$v_{(1^r,-1,0^s)}(x_{r+1}+1) > 0$,~$v_{(1^r,-1,0^s)}(x_j) > 0$ for~$1\leq j \leq s$
and~$v_{(1^r,-1,0^s)}(x_1-1) = 1$ (we normalize
this way in order to have formulas~(\ref{formules_valuations})).
We do the same for~$v_{(1^r)}$ and~$v_{(1^r,-1)}$.
With these notations, the Newton polygons of the preceding tabular
permit to prove that
every loop at~$1$ and the path form~$1$ to~$-1$ multiply the valuations by~$2$,
while the path from~$-1$ to~$0$ and the loop at~$0$ divides the valuations
by~$2$. More precisely, we have
\begin{equation}\label{formules_valuations}
v_{(1^r)}(x_r-1) = 2^{r-1},
\qquad
v_{(1^r,-1)}(x_{r+1}+1) = 2^{r},
\qquad
v_{(1^r,-1,0^s)}(x_{r+s+1}) = 2^{r-s}.
\end{equation}

Let us start with the point~$Q=(1^{n-2},-1,0)$, that is with~$r = n-2$ and~$s = 1$.
The point~$P = (1^{n-2},-1) \in X_{n-1}$ is smooth and thus corresponds to
a unique point of~$\widetilde{X}_{n-1}$, and the function~$x_1$ is an uniformizing parameter at~$P$. Let~$\OO_P$ be the local ring
at~$P$. There
is a unique singular point~$Q^\sharp\in X_n^\sharp$
above~$P \in X_{n-1}$, and above~$Q\in X_n$.
The local ring at~$Q^\sharp$
is~$\OO_{Q^\sharp}=\OO_P[x_n]$. Since~$v_P(x_{n-1}+1) = 2^{n-2}$, one has
$$
\left(\frac{x_n}{x_1^{2^{n-3}}}\right)^2
- 
\underbrace{\frac{x_{n-1}}{3x_{n-1}-1} \times \frac{x_{n-1}+1}{x_1^{2^{n-2}}}}_{\text{non zero at~$P$}}
= 0.
$$
The non-vanishing of the constant function shows that this equation is
smooth. This proves that~$\widetilde{\OO_P[x_n]}
= \OO_P\left[\frac{x_n}{x_1^{2^{n-3}}}\right]$. Thus
there are two points
in~$\widetilde{X_n}$ above~$Q^\sharp \in X_n^\sharp$,
or above~$(1^{n-2},-1,0) \in X_n$ and moreover~$\delta_{P^\sharp} = 2^{n-3}$.

We proceed by induction on~$s$ for the study of the singularity above the
general singular point~$P=(1^r,-1,0^s)\in X_n$ for~$r+s+1=n$.
We need to distinguish two cases.

If~$r\geq s\geq 1$, then
by induction, there are exactly~$2^{s-1}$ points of~$\widetilde{X}_{n-1}$
above~$(1^r,-1,0^{s-1})$. Let~$P$ one of them and
let~$\OO_P$ be the local ring at~$P$. The function~$x_1$
is still an uniformizing parameter at~$P$ and
above~$P$ there is only
one point~$Q^\sharp \in X_n^\sharp$, whose local ring is~$\OO_{Q^\sharp}=\OO_P[x_n]$.
Then~$v_P(x_{n-1}) = 2^{r-s+1}$, and the equation
$$
\left(\frac{x_n}{x_1^{2^{r-s}}}\right)^2
- 
\underbrace{\frac{x_{n-1}}{3x_{n-1}-1} \times \frac{x_{n-1}+1}{x_1^{2^{r-s+1}}}}_{\text{non zero at~$P$}}
= 0
$$
is smooth. This proves that there are
two points of~$\widetilde{X_n}$ above~$Q^\sharp \in X_n^\sharp$,
that is above~$(1^r,-1,0^s) \in X_n$, and that~$\delta_{P^\sharp} = 2^{r-s}$.

In the other case, that is if~$s > r \geq 1$,
then~$x_{n-1}$ turns to be an uniformizing element for all the
points of~$\widetilde{X}_{n-1}$ above~$P=(1^r,-1,0^{s-1})$, and
the corresponding points~$P^\sharp$ on~$X_n^\sharp$ are smooth
(and ramified over~$\widetilde{X}_{n-1}$).

In conclusion, we get by summation
$$
\Delta_n
=
2 \times \sum_{s=1}^{\lfloor\frac{n-1}{2}\rfloor} 2^{s-1} \times 2^{r-s}
=
\sum_{s=1}^{\lfloor\frac{n-1}{2}\rfloor} 2^{n-s-1},
$$
where the first factor $2$ comes from the other component of the singular
graph. The remaining of the computation is left to the reader.
\end{preuve}

\subsubsection{Lower bound for the number of points sequence}\label{LowerBound}

The last step is to prove that the tower has a large enough number of rational
points sequence over some~$\FF_{p^r}$, that is of 
size~$\sharp \widetilde{X}_n(\FF_{p^r}) = c\times 2^n + o(2^n)$ for some non-zero
constant~$c$ as requested in Proposition~\ref{Caract_Bonne_Tour}. To this end, we prove the following result.

\begin{prop}\label{number_of_points_sequence}
Let~$p\geq 5$ be a prime.
There exists~$r\geq 1$, and a finite set~$T \subset \PP^1(\FF_{p^r})$
with~$2(p-1)$ elements,
such that the graph~$\Graphe_T(\Cun,\Cdeux)$
is $2$-regular and the number of points sequence of the
tower~${\mathcal T}\left(\PP^1_{\FF_p}, y^2 = \frac{x^2+x}{3x-1}\right)$
satisfies
$$
\sharp X_n(\FF_q^r) \geq (p-1)\times 2^n.
$$
\end{prop}

By the correspondence
recalled in Proposition~\ref{PointsPaths}, if such a finite set~$T$
exists, then the
curve~$X_n$ contains at least~$2^n(p-1)$ points over~$\FF_{p^2}$.
The rest of this section is devoted to the proof of the existence
of this finite set~$T$, stated in~Lemma~\ref{Specialized_completeness_criterion}.
It is divided in few steps.

First, we experimentally compute this finite set
and its characteristic polynomial for few small primes~$p$. Second,
we observe that these polynomials
% for varying $p$ 
 for these small values of $p$ do lift to $\ZZ$, that is are reduction modulo~$p$ of some truncation of some integer
coefficients series. We enter this experimental integer sequence of first coefficients in the
database OEIS, which luckily returns a whole infinite integer sequence.
Finally, we prove that the reduction modulo~$p$ of some truncation of the
generating series of this infinite integer sequence do fulfill the functional equation required 
in Theorem~\ref{functional_completeness_criterion}. This implies
the regularness of the reciprocal image by~$f$ of the zeroes of the
truncations. The fact that the power series is
 closely related to a Gaussian hypergeometric function, hence satisfies some
second order linear differential equation, turns to be a crucial point here.

\medbreak

\begin{figure}
\centering
\begin{tikzpicture}[line width = 1pt,scale=1,baseline=(current bounding box.west)]
\fill (0,2) circle (3pt) node[left] {$\alpha^7$} ;
\fill (0,-2) circle (3pt) node[left] {$\alpha^3$} ;
\fill (2,1) circle (3pt) node[below left, xshift=-0.1cm] {$\alpha^{23}$} ;
\fill (2,-1) circle (3pt) node[above left, xshift=-0.1cm] {$\alpha^{21}$} ;
\fill (4,1) circle (3pt) node[below right, xshift=0.1cm] {$\alpha^{15}$} ;
\fill (4,-1) circle (3pt) node[above right, xshift=0.1cm] {$\alpha^{9}$} ;
\fill (6,2) circle (3pt) node[right] {$\alpha^{11}$} ;
\fill (6,-2) circle (3pt) node[right] {$\alpha^{19}$} ;
\draw[decoration={markings, mark=at position 0.5 with {\arrow{>}}},
      postaction={decorate}] (0,2) to [bend right] (2,1) ;
\draw[decoration={markings, mark=at position 0.5 with {\arrow{>}}},
      postaction={decorate}] (2,1) to [bend right] (2,-1) ;
\draw[decoration={markings, mark=at position 0.5 with {\arrow{>}}},
      postaction={decorate}] (2,-1) to [bend right] (0,-2) ;
\draw[decoration={markings, mark=at position 0.5 with {\arrow{>}}},
      postaction={decorate}] (0,-2) to [bend left] (0,2) ;
\draw[decoration={markings, mark=at position 0.5 with {\arrow{>}}},
      postaction={decorate}] (6,-2) to [bend right] (4,-1) ;
\draw[decoration={markings, mark=at position 0.5 with {\arrow{>}}},
      postaction={decorate}] (4,-1) to [bend right] (4,1) ;
\draw[decoration={markings, mark=at position 0.5 with {\arrow{>}}},
      postaction={decorate}] (4,1) to [bend right] (6,2) ;
\draw[decoration={markings, mark=at position 0.5 with {\arrow{>}}},
      postaction={decorate}] (6,2) to [bend left] (6,-2) ;
\draw[decoration={markings, mark=at position 0.5 with {\arrow{>}}},
      postaction={decorate}] (0,2) to [bend angle=12, bend right] (6,2) ;
\draw[decoration={markings, mark=at position 0.5 with {\arrow{>}}},
      postaction={decorate}] (6,2) to  [bend right] (0,2) ;
\draw[decoration={markings, mark=at position 0.5 with {\arrow{>}}},
      postaction={decorate}] (4,1) to [bend right] (2,1) ;
\draw[decoration={markings, mark=at position 0.5 with {\arrow{>}}},
      postaction={decorate}] (2,1) to [bend right] (4,-1) ;
\draw[decoration={markings, mark=at position 0.5 with {\arrow{>}}},
      postaction={decorate}] (4,-1) to [bend left] (0,-2) ;
\draw[decoration={markings, mark=at position 0.5 with {\arrow{>}}},
      postaction={decorate}] (0,-2) to [bend right] (6,-2) ;
\draw[decoration={markings, mark=at position 0.5 with {\arrow{>}}},
      postaction={decorate}] (6,-2) to [bend left] (2,-1) ;
\draw[decoration={markings, mark=at position 0.5 with {\arrow{>}}},
      postaction={decorate}] (2,-1) to [bend right] (4,1) ;
\draw (8,0.3) node[anchor = west] {$\FF_{25} = \FF_5(\alpha)$} ;
\draw (8,-0.3) node[anchor = west] {$\alpha^2 - \alpha + 2 = 0$} ;
\end{tikzpicture}
\caption{The $2$ regular complete component of~$\Graphe_\infty$ for~$p=5$}\label{figure_CC}
\end{figure}
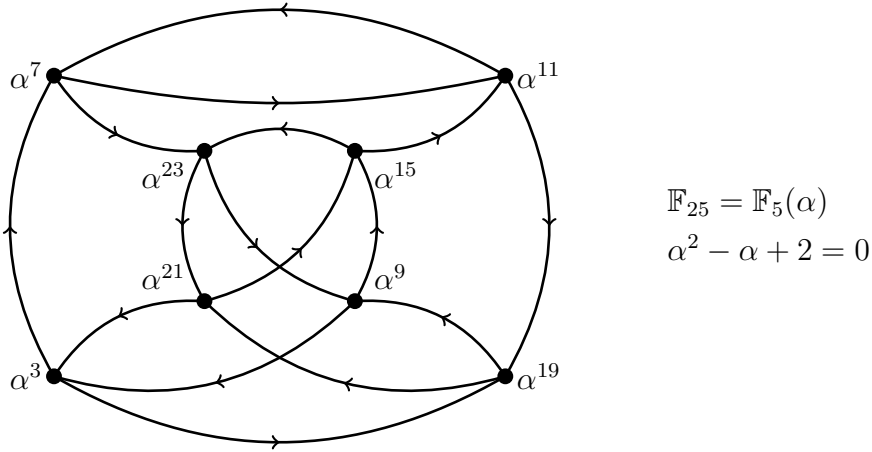

{\bfseries\noindent Experimental observation for few small primes $p$.~---}
It is an experimental observation, using {\tt magma} and {\tt sage},
that for small values of the prime~$p$,
the geometric graph~$\Graphe_\infty(\PP^1_{\FF_p}, y^2 = \frac{x^2+x}{3x-1})$
contains a finite $2$-regular component with~$2(p-1)$ vertices.
In figure~\ref{figure_CC}, we represent this component for~$p=5$.

\medbreak

{\bfseries\noindent Looking for a potentially splitting set for any prime $p$.~---}
Suppose that there exists, for any prime $p \geq 5$, some non-empty finite $2$-regular set $T=T_p$ in 
the geometric graph~$\Graphe_\infty(\PP^1_{\FF_p}, y^2 = \frac{x^2+x}{3x-1})$. Note that in this case, this set is
unique thanks to~\cite[Theorem~19]{HallouinPerretMMJ}.
Let~$\chi_p(x) \in \FF_p[x]$ be the characteristic
polynomial of the set of values of~$f$
(or~$g$) at the vertices of $T_p$
$$
\chi_p(x)
=
\prod_{P\in T_p} \left(x-x\left(f(P)\right)\right).
$$
It is an easy task using {\tt magma} to compute~$\chi_p(x)$
for small primes~$p$. Here is the table
for~$p \in \{5,7,11,13,17,23\}$:
$$
\begin{array}{l|l}
p & \chi_p(x) \\
\hline
5&
-1 + 2x + 2x^3 + x^4\\
\hline
7&
1 + 3x + x^2 + 2x^3 + 2x^4 - 2x^5 + x^6\\
\hline
11&
-1 - 3x - 4x^2 - 5x^3 - x^4 - 2x^6 - 2x^7 + x^8 + 4x^9 + x^{10}\\
\hline
13&
1 + 3x + 2x^2 + 2x^3 + 2x^4 - x^5 + 4x^6 + 4x^7 + 6x^8 + 2x^9 + 5x^{10} - 4x^{11} + x^{12}\\
\hline
17&
-1 - 3x + 2x^2 - 8x^3 + 7x^4 + 5x^5 + 4x^6 - 2x^7 + x^9 - x^{10} - 7x^{11} + 7x^{12} - 4x^{13}\\& - 8x^{14} + 6x^{15} + x^{16}\\
\hline
19&
1 + 3x - 4x^2 - 2x^3 - 7x^4 - 2x^5 - 5x^7 + 5x^8 + 7x^9 + 7x^{10} - 6x^{11} + 7x^{13} + 2x^{14}\\& - 3x^{15} + 3x^{16} - 6x^{17} + x^{18}\\
\hline
23&
-1 - 3x + 8x^2 - x^3 + 5x^4 - 7x^5 - 2x^6 - 9x^7 + 9x^8 - 9x^9 + 4x^{10} + 10x^{12} - 7x^{13}\\& - 6x^{14} + 8x^{15} - 7x^{16} - 2x^{17} + 3x^{18} - 10x^{19} + 7x^{20} + 8x^{21} + x^{22}
\end{array}
$$
We observe that the constant term is nothing but the Legendre
symbol~$\legendre{-3}{p}$. So all the
polynomials~$\legendre{-3}{p} \chi_p(x)$
have unitary constant coefficient.

These polynomials can be seen as the analogue for this tower, of the Deuring polynomials
for the tame optimal tower already touched on
in Subsection~\ref{TourGS}.
However, in contrast with this Subsection, we {\em a priori} do not know here any theoretical
polynomial which could play the part of the Deuring polynomials therein.
Fortunately, the Deuring
polynomials turn to be the modulo~$p$ truncations at degree~$(p-1)$ of a
power series in~$\ZZ[[x]]$.
One can expect that the same holds in our example.
More precisely, we want to prove that there
exists a power series~$H(x) = \sum_{n\geq 0} a_n x^n \in \ZZ[[x]]$, such that
for all primes~$p \geq 5$, one has
$$
\legendre{-3}{p} \chi_p(x)
=
\sum_{n= 0}^{p-1} a_n x^n
\bmod{p}.
$$
In other terms, we want to lift to~$\ZZ[[x]]$ the
polynomials~$\legendre{-3}{p} \chi_p(x) \in \FF_p[x]$. Thanks to the Chinese Remainder
Theorem and using more and more primes, 
one can lift
the experimentally known polynomials~$\legendre{-3}{p} \chi_p(x)$ for small $p$ modulo increasing integers.
For example using primes of the preceding tabular, one obtains
\begin{align*}
H(x)
=
1 &+ 3x + 15x^2 + 93x^3 + 639x^4 + 4653x^5 + 35169x^6 + 272835x^7 + 33065x^8\\
  &+ 322285x^9 + 438261x^{10} + 43884x^{11} + 40470x^{12} + 1755x^{13} - 2202x^{14}\\
  &+ 130x^{15} - 1327x^{16} - 44x^{17} + 20x^{18} + 10x^{19} - 7x^{20} - 8x^{21} - x^{22} + \cdots,
\end{align*}
each coefficient being respectively known modulo
\begin{align*}
&37182145,\, 37182145,\, 37182145,\, 37182145,\, 37182145,\, 7436429,\, 7436429,\, 1062347,\, 1062347,\\
&1062347,\, 1062347,\, 96577,\, 96577,\, 7429,\, 7429,\\
&7429,\, 7429,\, 437,\, 437,\, 23,\, 23,\, 23,\, 23.
\end{align*}
Requesting for the integer sequence~$1, 3, 15, 93, 639, 4\,653, 35\,169$
on the {\em Online Encyclopedia of Integer Sequences} \cite{OEIS},
we fortunately learn that these are the first terms of the sequence
$$
a_n = \sum_{k=0}^n \binom{n}{k}^2\binom{2k}{k}.
$$
The associated generating power series~$H(x)$ is moreover related to a Gaussian
hypergeometric function by
\begin{equation} \label{serie_H}
H(x)
=
\sum_{k=0}^\infty
\underbrace{\left(\sum_{k=0}^n \binom{n}{k}^2 \binom{2k}{k}\right)}_{a_n} x^n
=
\frac{1}{1-3x} \hypergeom{\frac{27x^2(1-x)}{(1-3x)^3}},
\end{equation}
where the  Gaussian hypergeometric part (see for instance~\cite[Chap. 5]{ConcrMath} for hypergeometric series) is
$$
\hypergeom{z}
=
\sum_{n=0}^\infty 
\frac{\left(\frac{1}{3}\right)_n \left(\frac{2}{3}\right)_n}
{\left(1\right)_n n!} z^n,
\qquad \hbox{so that}\qquad
\hypergeom{27z}
=
\sum_{n=0}^\infty 
\frac{(3n)!}{(n!)^3} z^n.
$$
Here, we denote as usual~$(a)_n \overset{\text{def}}{=} a(a+1) \cdots (a+n-1)$. We denote by~$H_p(x) \in \FF_p[x]$ the modulo~$p$ truncation in degree~$(p-1)$
of~$H$:
\begin{equation} \label{def_H_p}
H_p(x)
=
\sum_{n=0}^{p-1} a_n x^n \bmod p.
\end{equation}

\medbreak

{\noindent\bfseries Proof of the splitting behavior.~---} We have now reached the point where we do have, for any prime $p\geq 5$, a candidate for a $2$-regular
set, namely the set of roots of the explicit polynomial $H_p(x)$ defined in~\eqref{def_H_p}.
The encouraging point is that it is easily checked using {\tt magma} and {\tt sage} that for any tested value of $p$, this set
do corresponds to a $2$-regular component of the geometric graph.

\begin{lem}\label{Specialized_completeness_criterion}
Let~$p\geq 5$ be a prime and let~$T_0$ (depending on~$p$) be the set of roots
of the
polynomial~$H_p(x) \in \FF_p[x]$ defined in~(\ref{def_H_p}).
Then~$f^{-1}(T_0)$ is $2$-regular.
\end{lem}

\begin{preuve}
In order to apply the functional regularness criterion~Theorem~\ref{functional_completeness_criterion}
of which we keep notations, we choose as finite complete
set~$S$ the support of the singular graph
$$
S
=
\left\{0, \pm 1,  \pm\textstyle\frac{1}{3}, \infty\right\},
\qquad \hbox{so that}\qquad
S_0
=
\left\{0, 1, \textstyle\frac{1}{9}, \infty\right\}.
$$
The different divisors are
$$
\different_f(S_0) = \left[-\textstyle\frac{1}{3}\right] + \left[1\right]
\qquad\text{and}\qquad
\different_g(S_0) = \left[0\right] + \left[\infty\right],
$$
hence one can choose the function~$\rho$ to be
$$
\rho(x)
=
\frac{(x-1)\left(x+\frac{1}{3}\right)}{x} \in \FF_q(\PP^1).
$$
Since~$H_p(0) \equiv H_p(0) \equiv 1 \neq 0 \pmod{p}$,
the set~$T_0$ of roots of $H_p(x)$ in $\overline{\mathbb F_p}$ is disjoint from the ramification locus
$\{0, \infty\}$ of~$g$. 
Since~$\sharp T_0 = \deg(H_p) = (p-1)$, we
can try~$t=2$ and~$s=\frac{p-1}{2}$.
We put
$$
\varphi(x)
=
\frac{H_p(x)}{x(x-1)\left(x-\frac{1}{9}\right)}.
$$
After some easy computation similar to those in Subsection~\ref{TourGS}, we
see that the functional regularness criterion
specializes as follows in this example.
For the set~$f^{-1}(T_0)$ to be $2$-regular, it suffices that the
functional equation  
\begin{equation}\label{specialized_functional_equation}
\frac{1}{(3x-1)^{1-p}}
H_p
\left(
\frac{x^2+x}{3x-1}
\right)
\sim
H_p
\left(
x^2
\right)
\end{equation}
holds, which follows from Lemma~\ref{LemmeH} below.
\end{preuve}

\begin{lem}\label{LemmeCongr}
Let $n = \sum_{i=0}^r n_ip^i$ be the decomposition of $n$ in basis $p$. Then we have
$$a_n \equiv \prod_{i=1}^ra_{n_i} \pmod{p}.$$
\end{lem}

\begin{preuve} We recall Lucas formula~\cite[p. 271]{Dickson} for
binomial coefficients modulo a prime~$p$.
If~$n = n_0 + n_1 p + \cdots + n_r p^r$
and~$k = k_0 + k_1 p + \cdots + k_r p^r$ with~$p_i \in \{0, \dots, p-1\}$, then
$$
\binom{n}{k} \equiv \prod_{i=0}^r \binom{n_i}{k_i} \pmod{p}.
$$
From both computations
\begin{align*}
a_n & =  \sum_{k=\sum_{i=0}^r k_ip^i=0}^{n =\sum_{i=0}^r n_ip^i} \binom{\sum_{i=0}^r n_ip^i}{\sum_{i=0}^r k_ip^i}^2\binom{2k}{k}&\\
&\equiv \sum_{i=0}^r\sum_{k_i=0}^{p-1} \prod_{i=0}^r\binom{n_i}{k_i}^2\binom{2k}{k} \pmod{p} &\hbox{(by Lucas formula)}\\
\end{align*}
and
\begin{align*}
\prod_{i=1}^ra_{n_i} & \equiv   \sum_{i=0}^r\sum_{k_i=0}^{p-1} \prod_{i=0}^r\binom{n_i}{k_i}^2\prod_{i=0}^r\binom{2k_i}{k_i }\pmod{p}\\
\end{align*}
which takes into account that $\binom{n_r}{k_r} = 0$ if $k_r>n_r$, we have only to prove that for any $k=\sum_{i=0}^r k_ip^i$ with $0 \leq k_i\leq p-1$, equality
\begin{equation}\label{CongrAProuver}
\binom{2k}{k}\equiv \prod_{i=0}^r \binom{2k_i}{k_i} \pmod{p}
\end{equation}
holds. Suppose first that for any $i=0, \dots, r$, we have $0 \leq k_i \leq \frac{p-1}{2}$. Then $2k=\sum_{i=0}^r 2k_ip^i$ is the decomposition of $2k$ in basis $p$, and~\eqref{CongrAProuver} follows from Lucas formula. Suppose now that for at least one $i \in \{0, 1, \dots, r\}$, we have $\frac{p-1}{2}<k_i\leq p-1$. Then $k_i \leq p-1 < p \leq 2k_i$, hence as is well known $\binom{2k_i}{k_i} \equiv 0 \pmod{p}$. We now prove that in this case, $\binom{2k}{k} \equiv 0 \pmod{p}$, so that both sides in~\eqref{CongrAProuver} vanish.
 We use Legendre Theorem~\cite[p. 263]{Dickson} that the $p$-adic valuation of the factorial $n!$ of an integer $n$, written in basis $p$ as $n = \sum_{i=0}^rn_ip^i$, is given by
\begin{equation}\label{Legendre}
v_p(n!)=\frac{n-S_p(n)}{p-1},
\end{equation}
where $S_p(n) := \sum_{i=0}^r n_i$. The vanishing of $\binom{2k}{k}= \frac{(2k)!}{(k!)^2}$ modulo $p$ is equivalent to $v_p((2k)!) > 2v_p(k)$, hence by~\eqref{Legendre}, is equivalent to $S_p(2k)<2S_p(k)$. We observe that :
\begin{itemize}
\item for $i \in \{0, 1, \dots, r\}$ such that $0\leq k_i\leq \frac{p-1}{2}$, then $0\leq 2k_i \leq p-1$;
\item for $i \in \{0, 1, \dots, r\}$ such that $k_i=\frac{p-1}{2}+\ell_i$, with  $1\leq \ell_i \leq \frac{p-1}{2}$, then $2k_i= (2\ell_i-1)+p$, with $0\leq 2\ell_i-1 \leq p-1$.
\end{itemize}
Denote by $[N]_i$ the $i$-th digit of a composite integer $N$ in basis $p$, so that
\begin{equation}\label{v_p(fact)}
2S_p(k)-S_p(2k) = \sum_{i=0}^r 2k_i-[2k]_i.
\end{equation}
We deduce the following tabular, where contrib. means ``contribution to":

\medskip

\centerline{
\begin{tabular}{|c|clc|c|c|c|}
\hline
$k_i \in$ & $k_i$ & $2k_i$&  contrib. $[2k]_i$  & contrib. $[2k]_{i+1}$ & contrib. to~\eqref{v_p(fact)} \\
\hline 
$\{0, \dots, \frac{p-1}{2}\}$ & $k_i$ & $2k_i$&  $+2k_i$  & $+0$ & $+0$ \\
\hline 
$\{\frac{p-1}{2}+1, \dots, p-1\}$ & $\frac{p-1}{2}+\ell_i$ & $(2\ell_i-1)+p$&  $+(2\ell_i-1)$  & $+1$ & $+p-1>0$ \\ 
\hline
 \end{tabular}
}

\medskip
\noindent from which it follows that $v_p(\binom{2k}{k}) = 0$ if, and only if, there exists some $0\leq i \leq r$ such that $\frac{p-1}{2}<k_i$. The proof of Lemma~\ref{LemmeCongr} is complete.
\end{preuve}

We are now able to prove that the functional
equation~\eqref{specialized_functional_equation} holds. We use a Li's
trick \cite[\S7.2]{Li}, which
relies this truncated series modulo $p$ to the initial series,
and on the fact that hypergeometric functions are solutions of some second order linear differential equations.

\begin{lem} \label{LemmeH}
Let~$H(x) \in \ZZ[[x]]$ be the series defined in~(\ref{serie_H})
and for every prime~$p$, let~$H_p(x) \in \FF_p[x]$ be the degree~$(p-1)-th$
truncation of~$H(x)$ modulo~$p$ defined in~(\ref{def_H_p}).
\begin{enumerate}
\item\label{Li_Trick} The series~$H$ and its truncation~$H_p$ modulo~$p$ are related by
the relation
$$
H(x)^{1-p} \equiv H_p(x) \pmod{p}.
$$
\item\label{eqs_fonctionnelles} The series~$H$ and its truncation~$H_p$ satisfy the functional equations
$$
\frac{1}{1-3x}H\left(\frac{x^2+x}{3x-1}\right)
=
H(x^2)
\qquad\text{and}\qquad
\frac{1}{(1-3x)^{1-p}}H_p\left(\frac{x^2+x}{3x-1}\right)
=
H_p(x^2).
$$
\end{enumerate}
\end{lem}

\begin{preuve} Point~$\ref{Li_Trick}$.
We have the following congruences modulo~$p$:
\begin{align*}
H_p(x) \times H_p(x)^p \times H_p(x)^{p^2}\times \dots
&\equiv H_p(x) \times H_p(x^p) \times H_p(x^{p^2})\times \dots\\
& \equiv \prod_{i=0}^{\infty} \left(\sum_{n_i=0}^{p-1}a_{n_i}x^{n_ip^i}\right)\\
&\equiv \sum_{r \in \NN ; 0\leq n_0, \dots, n_r \leq p-1}\prod_{i=1}^ra_{n_i} x^{n_0+pn_1+ \cdots+ n_rp^r} \\
&\equiv \sum_{n \in \NN}a_{n} x^{n} \qquad \hbox{ (by Lemma~\ref{LemmeCongr})} \\
& \equiv H(x),
\end{align*}
where the first product converges to an invertible function in~$\ZZ_p[[x]]$. It follows that
$$H(x)^{1-p} = \frac{H(x)}{H(x)^p}  \equiv H_p(x) \pmod{p},$$
which proves~$\ref{Li_Trick}$. Now, Gaussian hypergeometric functions are known
to be solution of some second order linear differential equation. As for the hypergeometric
geometric function~$F(x) = \hypergeom{x}$, it satisfies the equation
$$
x(1-x)F''(x)
+
(1-2x)F'(x)
-
\frac{2}{9}F(x)
=0
$$
(see~\cite[Ex. 5.108 p. 221]{ConcrMath}). We then deduce (the details of the computations are left to the reader) two
second order linear differential equations respectively satisfied by the
functions~$x \mapsto \frac{1}{1-3x}H\left(\frac{x^2+x}{3x-1}\right)$
and~$x \mapsto H(x^2)$.
These two equations turn to be proportional.
Since the two preceding functions have same value and derivative at zero,
they must be equal. This complete the proof of the first functional
equation. To prove the second one, it suffices to raise the first one to the power~$(p-1)$
and to use point~$\ref{Li_Trick}$. This completes the proofs of~$\ref{eqs_fonctionnelles}$.
\end{preuve}

\subsubsection{The last question}
%We have not gave an answer to the following important question of the (unique thanks to Theorem of~\cite{HallouinPerretMMJ})
%value of $r$, such that
%\begin{quote}
%{\bfseries Question ---} For which~$r\geq 1$ do we
%have~$\lambda_r\left({\mathcal T}\left(\PP^1, y^2 = \frac{x^2+x}{3x-1}\right)\right) > 0$?
%\end{quote} 
%In fact we do have the response to this question. This is~$r = 2$, but we
%unfortunately need to remark that our tower is {\em not new} to prove this fact!

We have not yet answered the important question of the
value of $r$, such that $\lambda_r\left({\mathcal T}\left(\PP^1, y^2 = \frac{x^2+x}{3x-1}\right)\right) > 0$. The {\tt magma} experiments for small values of $p\geq 5$ show that, at least for these values of $p$, one has $r=2$. Unfortunately, we were not able to prove this. But a close look at Elkies article~\cite{ElkiesET} leads to the following Proposition, showing that this tower is actually not new.

\begin{prop}\label{Sur_Fp2}
The recursive tower~${\mathcal T}\left(\PP^1, y^2 = \frac{x^2+x}{3x-1}\right)$
is isomorphic to the modular tower $(X_0(3\cdot 2^n))_{n\geq 2}$ described
by Elkies \cite{ElkiesET}. This tower is asymptotically good, and even optimal,
over~$\FF_{p^2}$ for every prime~$p$ outside~$\{2,3\}$.
\end{prop}

\begin{preuve}
The model of the tower~$(X_0(3\cdot 2^n))_{n\geq 2}$ given by
Elkies \cite[formula~(45)]{ElkiesET} is the recursive tower with base
curve~$\Cun = \PP^1$ and correspondence~$\Cdeux_{f_{E},g_{E}}$
defined by the two functions~$f_{E}(x) = x^2$ and~$g_{E}(y) = \frac{y^2+3y}{y-1}$.
One easily verifies that~$f = \sigma \circ f_{E} \circ \tau$
and~$g = \sigma \circ g_{E} \circ \tau$
for~$\tau(x) = \frac{3x+1}{x-1}$ and~$\sigma(x) = \frac{x-1}{x-9}$.
\end{preuve}

%------------------------------------
%\bibliographystyle{alpha}
\bibliographystyle{amsalpha}
%\bibliography{Lenstra}
\providecommand{\bysame}{\leavevmode\hbox to3em{\hrulefill}\thinspace}
\providecommand{\MR}{\relax\ifhmode\unskip\space\fi MR }
\providecommand{\MRhref}[2]{%
  \href{http://www.ams.org/mathscinet-getitem?mr=#1}{#2}
}
\providecommand{\href}[2]{#2}

\bigbreak

\noindent
\begin{minipage}[t]{0.6\textwidth}
Hallouin Emmanuel ({\tt hallouin@univ-tlse2.fr})

Perret Marc ({\tt perret@univ-tlse2.fr})
\end{minipage}
\hfill
\begin{minipage}[t]{0.4\textwidth}
Universit\'e Toulouse Jean Jaur\`es

5, all\'ees Antonio Machado

31058 Toulouse cedex

France
\end{minipage}

\end{document}